\documentclass[11pt,twoside]{article}
\usepackage{enumerate}
\usepackage{amsmath}
\usepackage{amssymb}
\usepackage{amsthm}
\usepackage{amscd}
\DeclareSymbolFontAlphabet{\Bbb}{AMSb}
\newcommand{\id}{\text{{\rm id}}}

\newcommand{\norm}[1]{ \|#1 \| }
\newcommand{\bignorm}[1]{\left \|#1 \right \| }
\newcommand{\convex}[1]{\mathbf{M^{(#1)}}}
\newcommand{\concave}[1]{\mathbf{M_{(#1)}}}
\newcommand{\cotype}[1]{\mathbf{C_{#1}}}
\newcommand{\type}[1]{\mathbf{T_{#1}}}

\newcommand{\K}{\Bbb{K}}
\newcommand{\C}{\Bbb{C}}

\newcommand{\LL}{{\cal L}}
\newcommand{\MM}{{\cal M}}
\newcommand{\PP}{\Pi}
\newcommand{\I}{{\cal I}}

\newcommand{\Id}{\hookrightarrow}

\newcommand{\ui}{\mathcal{S}}
\theoremstyle{definition}
\newtheorem{defin}{Definition}[section]

\newtheorem{ex}[defin]{Example}

\theoremstyle{plain}
\newtheorem{lemma}[defin]{Lemma}
\newtheorem{theo}[defin]{Theorem}
\newtheorem{cor}[defin]{Corollary}
\newtheorem{prop}[defin]{Proposition}

\theoremstyle{remark}
\numberwithin{equation}{section}
\begin{document}
\title{\bf Complex interpolation of spaces of operators on 
$\boldsymbol{\ell_1}$} 
\author{Andreas Defant and Carsten
Michels 
\\ \small\sl Carl von Ossietzky Universit\"at Oldenburg,
Fachbereich  Mathematik 
\\ \small\sl Postfach 2503, D-26111
Oldenburg, Germany 
\\ \small\tt defant\small\rm @\small\tt
mathematik.uni-oldenburg.de 
\\ \small\tt michels\small\rm
@\small\tt mathematik.uni-oldenburg.de } 
\date{} 
\maketitle
\begin{abstract} 
Within the theory of complex interpolation and
$\theta$-Hilbert spaces we extend classical results of
Kwapie\'{n}  on absolutely  $(r,1)$-summing operators on
$\ell_1$ with values in $\ell_p$ as well as their natural
extensions for mixing operators invented by Maurey. Furthermore,
we show that for  $1<p<2$ every operator $T$ on $\ell_1$ with
values in $\theta$-type~$2$  spaces, $\theta=2/{p'}$, is
Rademacher~$p$-summing. This is another  extension of
Kwapie\'{n}'s results, and by an extrapolation procedure  a
natural supplement to a statement of Pisier. 
\end{abstract} 
\setcounter{section}{-1} 
\section{Introduction} 
Kwapie\'{n} in
\cite{kwapien} showed that for $1 \le p \le \infty$ and  $1 \le
r \le 2$ defined by \mbox{$1/r=1-|1/2-1/p|$} every continuous 
and linear operator  on $\ell_1$ with values in $\ell_p$ is
$(r,1)$-summing, i.\,e., maps unconditionally summable into
absolutely $r$-summable sequences, and  Pisier in
\cite{pisier79}  proved that this result also holds whenever
$\ell_p$ ($1 \le p \le 2$) is replaced by an arbitrary
$p$-convex and $p'$-concave Banach function space. Carl and the
first author in  \cite{CD92} gave a generalization of
Kwapie\'{n}'s result  within the  framework of mixing operators:
For $2 \le s \le \infty$ such that  $1/s=|1/2-1/p|$ every
operator $T:\ell_1 \rightarrow \ell_p$ is  $(s,1)$-mixing,
i.\,e. the composition of $T$ with an arbitrary $s$-summing
operator  on $\ell_p$ is $1$-summing.  While Kwapie\'{n} and
Pisier used interpolation techniques  (e.\,g. the
Three-Lines-Theorem together with results of  Orlicz ($p=1$ or
$p=\infty$) and Grothendieck ($p=2$)), Carl and Defant used a 
certain tensor product trick.  
\par In this paper we suggest a
systematic approach to all these results within  the framework
of complex interpolation and $\theta$-Hilbert spaces which  for
example allows to replace $\ell_p$ in Kwapie\'{n}'s result by
the  Schatten class $\ui_p$ and also covers the well-known
results of Mitiagin \cite{kwapien} and  \cite{CD92} on the
coincidence of summing/mixing operators on $\ell_2$ and Schatten
classes. Furthermore,  we show that for $1 \le p  \le 2$ every
operator on $\ell_1$ with values in a $p$-convex Banach function
space  $X$ with finite concavity is Rademacher $p$-summing; in
other terms, each  operator $T:\ell_1 \rightarrow X$ maps weakly
$p$-summable sequences into almost  unconditionally convergent
sequences. This result follows by extra- and  interpolation
between the two border cases $p=1$ (trivial case) and  $p=2$
($X$ has type~2 if and only if $X$ is $2$-concave and has
non-trivial  concavity if and only if every operator $T: \ell_1
\rightarrow X$ is  Rademacher $2$-summing). Moreover, it is a
natural supplement to the  above mentioned result of
Pisier---note that in contrast to his result our  result  on
Rademacher $p$-summing operators does not depend on the exact 
degree of concavity of the image space.  
\par Finally, we  refer
the reader to a recent paper by  Cerd\`{a} and Masty\l o
\cite{cerda} where an extension of Kwapie\'{n}'s  result  within
the framework of Calder\'{o}n--Lozanovski\u{\i} spaces, in
particular Orlicz spaces, is given. 
\par We  use standard
notation and notions from Banach space theory, as presented
e.\,g. in  \cite{lt77} and \cite{lt}.   For $1 \le p \le \infty$
the number $p'$ is defined by $1/p + 1/{p'} =1$. If $E$ is a
Banach space, then $B_E$ is its (closed) unit ball and $E'$ its
dual.  As usual $\LL(E,F)$ denotes the Banach space of all
(bounded and linear) operators from $E$ into $F$ endowed with
the operator norm $\norm{\cdot}$.  For $1 \le p \le 2 \le q <
\infty$ we denote by $\type{p}(E)$ and  $\cotype{q}(E)$ the
(Rademacher) type~$p$ constant and the  (Rademacher) cotype~$q$
constant of a Banach space with these properties. 
\par  Let
$(\Omega,\Sigma,\mu)$ be a $\sigma$-finite and complete measure
space, and denote all $\mu$-a.e. equivalence classes of
real-valued measurable functions on $\Omega$ by $L_0(\mu)$. A
Banach space $X=X(\mu)$ of  (equivalence classes of) functions
in $L_0(\mu)$ is said to be a  Banach function space if it
satisfies the following conditions: 
\\[10pt] (I) If $|f| \le
|g|$, with $f \in L_0(\mu)$ and $g \in X(\mu)$, then $f \in
X(\mu)$ and $\norm{f}_X \le \norm{g}_X$. 
\\[5pt] (II) For every
$A \in \Sigma$ with $\mu(A)<\infty$ the characteristic  function
$\chi_A$ of $A$ belongs to $X(\mu)$. 
\\[10pt] For $1 \le p \le q
\le \infty$ we denote by $\convex{p}(X)$ and  $\concave{q}(X)$
the $p$-convexity constant and the $q$-concavity constant of a
Banach function space $X$ with these properties  (see e.\,g.
\cite{lt}), and with $X(\C)$ its natural complexification. 
\par
For all information on Banach operator ideals see e.\,g.
\cite{df},  \cite{djt}, \cite{pietsch} and \cite{tj}, and for
the theory of tensor norms on  tensor products of Banach spaces
\cite{df}. If $E$ is a symmetric Banach  sequence space, then we
denote by $\ui_E$ its associated unitary ideal,  i.\,e. the
space of all compact $T: \ell_2 \rightarrow \ell_2$ for which  
the sequence $(s_n(T))_n$ of singular numbers belongs to $E$,
equipped  with the norm $\norm{T}_{\ui_E}:=
\norm{(s_n(T))_n}_E$.  
\par  
We point out that, since we
extensively use complex  interpolation, the underlying field
$\K$ is always $\C$---important  exceptions are mentioned
explicitly. However, many of our main results can be  easily
transferred to the real case; we leave this to the  reader. 
\section{Complex interpolation and approximation} 
For an
introduction to complex interpolation theory see \cite{BL} and
\cite{kps}. If we speak of a finite-dimensional interpolation
couple $[E_0,E_1]$, we always  assume that $E_0$ and $E_1$ have
the same finite dimension. For $0<\theta<1$ and a complex
interpolation couple we denote by $[E_0,E_1]_\theta$ the 
associated complex interpolation space, and if $E_1$ is a
Hilbert space,  then $[E_0,E_1]_\theta$ is called a
$\theta$-Hilbert space; this notion is  due to Pisier. Note that
one can always assume that $E_0 \cap E_1$ is dense in $E_0$ and
$E_1$. 
\par Complex interpolation of finite-dimensional Banach
spaces is rather more  interpolation of norms than interpolation
of spaces which makes it sometimes easier to explain an idea,
without struggling with topological questions.  Therefore we try
to prove most of our results in the finite-dimensional case
first (with ``uniform'' control of certain constants), and then
we lift them to the infinite-dimensional case with the help of
an  ``approximation'' lemma which is proved in this first
section. 
\\[10pt]   For this purpose we introduce the notion of
a ``cofinal  interpolation triple'': If $[E_0,E_1]$ is an
interpolation couple,  $E \subset E_0 \cap E_1$ a subspace which
is dense in $E_0,E_1$, and  $\mathcal{B} \subset FIN(E)$ (where
$FIN(E)$ denotes the collection of all finite-dimensional
subspaces of $E$) is cofinal (i.\,e. for every  $G \in FIN(E)$
there exists $M \in \mathcal{B}$ with $G \subset M$), then the
triple $([E_0,E_1],E,\mathcal{B})$ is called a {\em cofinal 
interpolation triple.} For $M \in FIN(E)$ we denote by $M_0$ and
$M_1$ the subspace $M$ of $E_0$ and $E_1$ endowed with the
induced norm,  respectively. The following lemma seems to be
folklore in Russian literature; a proof can be found in e.\,g.
\cite[4.1]{kouba} or \cite[1.3]{michdiss}. 
\begin{lemma}
\label{lift} 
Let $([E_0,E_1],E,\mathcal{B})$ be a cofinal
interpolation triple and  $0<\theta<1$. Then for each
$\varepsilon>0$ and $G \in FIN(E)$ there  exists $M \in
\mathcal{B}$ such that $G \subset M$ and for all $x \in G$ 
$$
(1-\varepsilon) \cdot \norm{x}_{[M_0,M_1]_\theta} \le 
\norm{x}_{[E_0,E_1]_\theta} \le \norm{x}_{[M_0,M_1]_\theta}. 
$$
\end{lemma} 
Next we prove the approximation lemma announced
above. 
\begin{samepage}   
\begin{lemma} 
\label{l1lift} 
Let
$([F_0,F_1],F,\mathcal{B})$ be a cofinal interpolation triple,
$0<\theta<1$  and $(\mathcal{A},A)$ a maximal Banach operator
ideal with $$ c_\theta:=\sup_n \sup_{M\in \mathcal{B}}
\norm{\LL(\ell_1^n, [M_0,M_1]_\theta) \Id
\mathcal{A}(\ell_1^n,[M_0,M_1]_\theta)} <\infty. $$ Then 
$$
\LL(\ell_1,[F_0,F_1]_\theta)=\mathcal{A}(\ell_1,[F_0,F_1]_\theta). 
$$ 
\end{lemma} 
\end{samepage} 
\proof Denote by  $\varepsilon$
the injective tensor norm, by  $\alpha$ the finitely generated
tensor norm associated to  $(\mathcal{A},A)$ and by
$\overset{\leftarrow}{\alpha}$  its cofinite hull in the sense
of \cite[17.3 and 12.4]{df}. We prove that 
\begin{equation}
\label{cofinite} 
\norm{\ell_\infty^n \otimes_\varepsilon
(F,\norm{\cdot}_{[F_0,F_1]_\theta}) \Id \ell_\infty^n
\otimes_{\overset{\leftarrow}{\alpha}} 
(F,\norm{\cdot}_{[F_0,F_1]_\theta})}\le c_\theta. 
\end{equation}
Then by density (see \cite[13.4]{df}) 
$$ \norm{\ell_\infty^n
\otimes_\varepsilon [F_0,F_1]_\theta \Id \ell_\infty^n
\otimes_{\overset{\leftarrow}{\alpha}}  [F_0,F_1]_\theta}\le
c_\theta, 
$$ 
hence the claim follows by the Embedding Theorem
\cite[17.6]{df} and  local techniques \cite[23.1]{df}. \par For
$z \in \ell_\infty^n \otimes F$ choose by Lemma~\ref{lift} a
subspace $M \in \mathcal{B}$  such that  $z \in \ell_\infty^n
\otimes M$ and  \mbox{$\norm{(M,\norm{\cdot}_{[F_0,F_1]_\theta}
) \Id [M_0,M_1]_\theta } \le 1+\varepsilon$}. Then by  the
mapping properties of $\overset{\leftarrow}{\alpha}$ and 
$\varepsilon$ and the fact that  the injective tensor norm
respects subspaces we have 
$$ \norm{z}_{\ell_\infty^n
\otimes_{\overset{\leftarrow}{\alpha}}
(F,\norm{\cdot}_{[F_0,F_1]_\theta})} \le 
\norm{z}_{\ell_\infty^n \otimes_{\overset{\leftarrow}{\alpha}} 
[M_0,M_1]_\theta} \quad \text{and} \quad \norm{z}_{\ell_\infty^n
\otimes_\varepsilon  [M_0,M_1]_\theta}  \le (1+\varepsilon)
\cdot \norm{z}_{\ell_\infty^n \otimes_\varepsilon 
(F,\norm{\cdot}_{[F_0,F_1]_\theta})}, 
$$ 
hence \eqref{cofinite}
follows from the assumption and the  Embedding Theorem. \qed
\par In our applications the approximation lemma will be
combined with  what we call the ``interpolation trick'', due to
Kwapie\'{n} and based on complex  interpolation of vector-valued
$\ell_p$'s \cite[5.1.2]{BL}:  
$$
\LL(\ell_1^n,[M_0,M_1]_\theta)=\ell_\infty^n([M_0,M_1]_\theta) 
=[\ell_\infty^n(M_0),\ell_\infty^n(M_1)]_\theta =
[\LL(\ell_1^n,M_0), \LL(\ell_1^n,M_1)]_\theta 
$$
(isometrically).  
\section[Summing and mixing operators on
$\ell_1$]{Summing and mixing  operators on
$\boldsymbol{\ell_1}$} 
\par For all information on  summing
and mixing operators see e.\,g. \cite{pietsch}, \cite{djt} and
\cite{df}. An operator $T \in \LL(E,F)$ is called absolutely
$(r,p)$-summing  $(1 \le p \le r \le \infty)$ if there is a
constant $\rho \ge 0$ such that 
$$ 
\bigl ( \sum_{i=1}^n
\norm{Tx_i}^r \bigr )^{1/r} \le \rho \cdot \sup \biggl \{ \bigl
( \sum_{i=1}^n |\langle x',x_i \rangle |^p  \bigr )^{1/p} \, |
\, x' \in B_{E'} \biggr \} 
$$ 
for all finite sets of elements
$x_1, \dots, x_n \in E$ (with the obvious modifications for $p$
or $r$ $=\infty$). In this case, the infimum over all possible
$\rho \ge 0$ is denoted by $\pi_{r,p}(T)$,  and the maximal
Banach operator ideal of all absolutely $(r,p)$-summing 
operators by $(\PP_{r,p},\pi_{r,p})$; the special case $r=p$
gives the ideal $(\PP_p, \pi_p)$ of all absolutely $p$-summing
operators. 
\par An operator $T \in \LL(E,F)$ is called
$(s,p)$-mixing  $(1 \le p \le s \le \infty)$ whenever its
composition with an  arbitrary operator $S \in \PP_s(F,Y)$ is
absolutely $p$-summing; with the norm 
$$ 
\mu_{s,p}(T) :=
\sup\{\pi_p(ST) \, | \, \pi_s(S) \le 1 \} 
$$ 
the class
$\MM_{s,p}$ of all $(s,p)$-mixing operators forms again a
maximal Banach operator ideal. Obviously, $(\MM_{p,p},
\mu_{p,p}) = (\LL, \norm{\cdot})$ and $(\MM_{\infty,p},
\mu_{\infty,p}) = (\PP_p, \pi_p)$.   Recall that due to
\cite{maurey} (see also \cite[32.10--11]{df}) summing  and
mixing operators are closely related:  
\begin{equation}
\label{inclusion} 
(\MM_{s,p}, \mu_{s,p}) \subset (\PP_{r,p},
\pi_{r,p}) \qquad \text{for }    1/s + 1/r = 1/p, 
\end{equation} 
and ``conversely''  
\begin{equation}
\label{reverse} (\PP_{r,p}, \pi_{r,p}) \subset (\MM_{s_0,p},
\mu_{s_0,p}) \qquad \text{for } 1 \le p \le s_0 < s \le \infty
\text{ and } 1/s + 1/r = 1/p.   
\end{equation} 
Moreover, it is
known that the identity map $\id_X$ of a cotype~$2$ space  $X$
is  $(2,1)$-mixing and therefore every $(s,2)$-mixing operator
on a cotype~2 space is even $(s,1)$-mixing (see again
\cite{maurey} and \cite[32.2]{df}). Finally, a quick 
investigation of \cite{tj70} shows that if $1 \le q \le 2 \le r
\le \infty$ with $1/r=1/q-1/2$ are given, then
$\Pi_{r,2}(X,\cdot)=\Pi_{q,1}(X,\cdot)$  for every Banach space
$X$ such that  $\id_X$ is  $(2,1)$-mixing; by the above, this
holds in particular for cotype~2 spaces, hence, most of our main
results in this article can also be formulated  in terms of
$(q,1)$-summing/mixing norms instead of $(q,2)$-summing/mixing
norms. 
\par The following theorem is an extension of the results by Kwapie\'{n} and
\cite{CD92}---it is an unpublished result from \cite{referee} and 
Lermer \cite{lermer} (however with proofs different from the one presented
 in this article).
\begin{theo} 
\label{theta1}
Let $F$ be a $\theta$-Hilbert space, $0 < \theta < 1$. Then  
$$
\LL(\ell_1, F) = \MM_{\frac{2}{1-\theta},2}(\ell_1,F) = 
\PP_{\frac{2}{\theta},2}(\ell_1,F). 
$$ 
\em Since every
$2$-summing operator on $\ell_1$ is $1$-summing this implies 
\begin{equation} 
\label{thetaa} 
\LL(\ell_1, F) =
\MM_{\frac{2}{1-\theta},1}(\ell_1,F) = 
\PP_{\frac{2}{1+\theta},1}(\ell_1,F). 
\end{equation} 
\end{theo}
The proof is based on the following interpolation theorem for spaces of 
mixing
operators in combination with Lemma~\ref{l1lift}.  
\begin{prop}
\label{mixing2} 
Let $[E_0, E_1]$ and $[F_0, F_1]$ be
finite-dimensional interpolation  couples   and  $2 \le s_0, s_1
\le \infty$. Then for $0\le \theta \le 1$ and $2 \le s_\theta
\le \infty$ defined by \mbox{$1/{s_\theta} = (1-\theta)/{s_0} +
\theta/{s_1}$} 
$$\norm{[{\cal M}_{s_0,2}(E_0,F_0), {\cal
M}_{s_1,2}(E_1,F_1)]_\theta \Id  {\cal
M}_{s_\theta,2}([E_0,E_1]_\theta,[F_0,F_1]_\theta)}  \le
d_\theta[E_0,E_1], 
$$ 
where $d_\theta[E_0,E_1]:=\sup_m
\norm{\LL(\ell_2^m,[E_0,E_1]_\theta)  \Id
[\LL(\ell_2^m,E_0),\LL(\ell_2^m,E_1)]_\theta}$. 
\\[10pt]  {\em
Note that in the case $E_0=E_1$ (isometrically), one  trivially
has $d_\theta[E_0,E_1]=1$.} 
\end{prop} 
\proof Consider for
$\eta=0,1$  the  mappings 
$$ 
\begin{array}{lccccccc}
\Phi_\eta^{n,m}: &\MM_{s_\eta,2}(E_\eta, F_\eta) &\times
&\ell_{s_\eta}^n (F_\eta') & \times &\LL(\ell_2^m,E_\eta)
&\longrightarrow & \ell_2^m(\ell_{s_\eta}^n) \\ & T & \times
&(y_1', \dots , y_n') &\times & S &\longmapsto & ((\langle
y_k',TSe_j \rangle )_k)_j 
\end{array}, 
$$ 
where $(e_j)$ denotes
the canonical basis in $\C^m$. By the discrete characterization
of the mixing norm (see \cite{maurey} or \cite[32.4]{df})
$\mu_{s_\eta,2}(T:E_\eta \rightarrow F_\eta)$ is the infimum
over all $c \ge 0$ such that for all $n,m$, all  $y_1',
\ldots,y_n' \in F_\eta'$ and all $x_1,\ldots,x_m \in E_\eta$ 
$$
\left (\sum_{j=1}^m \left (\sum_{k=1}^n |\langle y_k',Tx_j
\rangle|^{s_\eta} \right )^{2/{s_\eta}} \right )^{1/2} \le c
\cdot \left ( \sum_{k=1}^n  \norm{y_k'}_{F_\eta'}^{s_\eta}
\right )^{1/{s_\eta}} \cdot  \sup_{x' \in B_{E_\eta'}} \left (
\sum_{j=1}^m |\langle x',x_j \rangle|^2  \right )^{1/2}. 
$$
Since for each $S= \sum_{j=1}^m e_j \otimes x_j \in
\LL(\ell_2^m,E_\eta)$ $$ \norm{S} =\sup_{x' \in B_{E_\eta'}}
\left ( \sum_{j=1}^m |\langle x',x_j \rangle|^2  \right )^{1/2},
 $$ it clearly follows that $\norm{\Phi_\eta^{n,m}}\le 1$. Then
for the interpolated mapping  
\begin{multline*}
[\Phi_0^{n,m},\Phi_1^{n,m}]_\theta : [\MM_{s_0,2}(E_0,F_0),
\MM_{s_1,2}(E_1,F_1)]_\theta \times [\ell_{s_0}^n(F_0'),
\ell_{s_1}^n(F_1')]_\theta  \\ \times
[\LL(\ell_2^m,E_0),\LL(\ell_2^m,E_1)]_\theta  \rightarrow
[\ell_2^m(\ell_{s_0}^n),\ell_2^m(\ell_{s_1}^n)]_\theta,
\end{multline*} 
by multilinear interpolation (see e.\,g.
\cite[4.4.1]{BL}) we also have  $\norm{[\Phi_0^{n,m},
\Phi_1^{n,m}]_\theta} \le 1$. It follows that for each $T:
[E_0,E_1]_\theta \rightarrow [F_0,F_1]_\theta$, each $S \in 
\LL(\ell_2^m,[E_0,E_1]_\theta)$ and $y_1', \ldots, y_n' \in 
[F_0,F_1]_\theta'$ 
\begin{multline*} 
\left ( \sum_{j=1}^m \left
( \sum_{k=1}^n | \langle y_k',TSe_j \rangle  |^{s_\theta} \right
)^{2/{s_\theta}} \right )^{1/2} 
\\ \begin{split} & \le
\norm{T}_{[\MM_{s_0,2}(E_0,F_0),\MM_{s_1,2}(E_1,F_1)]_\theta} 
\cdot \norm{S}_{[\LL(\ell_2^m,E_0),\LL(\ell_2^m,E_1)]_\theta}
\cdot
\norm{(y_k)_k}_{[\ell_{s_0}^n(F_0'),\ell_{s_1}^n(F_1')]_\theta}
\\ & \le d_\theta[E_0,E_1] \cdot
\norm{T}_{[\MM_{s_0,2}(E_0,F_0), \MM_{s_1,2}(E_1,F_1)]_\theta}
\cdot \norm{S}_{\LL(\ell_2^m, [E_0,E_1]_\theta)} \cdot
\norm{(y_k)_k}_{\ell_{s_\theta}^n ([F_0,F_1]_\theta')},
\end{split} 
\end{multline*} 
hence 
\begin{samepage}  
$$
\norm{T}_{\MM_{s_\theta,2}([E_0,E_1]_\theta,[F_0,F_1]_\theta)}
\le  d_\theta[E_0,E_1] \cdot
\norm{T}_{[\MM_{s_0,2}(E_0,F_0),\MM_{s_1,2} (E_1,F_1)]_\theta}.
$$  
\qed 
\end{samepage} 
\par  {\em Proof} of
Theorem~\ref{theta1}: Let $F=[F_0,F_1]_\theta$ where $F_1$ is a
Hilbert space and $F_0 \cap  F_1$ is dense in $F_0$ and $F_1$
(this implies that $([F_0,F_1], F_0 \cap F_1,  FIN(F_0 \cap
F_1))$ is a cofinal interpolation triple), and let   $M \in
FIN(F_0 \cap F_1)$. Then  by the Little Grothendieck Theorem
(see e.\,g. \cite[11.11]{df}) 
$$ \norm{\LL(\ell_1^n,M_1) \Id
\MM_{\infty,2}(\ell_1^n,M_1)=\Pi_2(\ell_1^n, M_1)} \le K_{LG},
$$ 
and trivially $\norm{\LL(\ell_1^n,M_0) \Id
\MM_{2,2}(\ell_1^n,M_0)  =\LL(\ell_1^n,M_0)} \le 1$, hence, by
the usual interpolation theorem together with
Proposition~\ref{mixing2} and Kwapie\'{n}'s  interpolation
trick, we have 
$$ \norm{\LL(\ell_1^n,[M_0,M_1]_\theta
)=[\LL(\ell_1^n,M_0),\LL(\ell_1^n,M_1)]_\theta \Id
\MM_{2/(1-\theta),2} (\ell_1^n,[M_0,M_1]_\theta)} \le
K_{LG}^\theta 
$$ 
($d_\theta[\ell_1^n,\ell_1^n]=1$). The claim
now follows by  Lemma~\ref{l1lift}. \qed 
\par We illustrate the
preceding theorem by several examples and a corollary (see the  remarks after
the proof of the corollary for credits): 
\begin{ex}
\label{example1a} 
For $1<p<\infty$, $p \neq 2$, and $\theta = 1-
|1-2/p|$ it is known that $\ell_p$ and $\ui_p$ are
$\theta$-Hilbert  spaces (see e.\,g. \cite[5.1.1]{BL} and
\cite[Satz~8]{pt} together with the complex reiteration theorem
\cite[4.6.1]{BL}), hence we obtain  
\begin{equation}
\label{lpkwapien}
\LL(\ell_1,\ell_p)=\MM_{2/{|1-2/p|},2}(\ell_1,\ell_p)=\Pi_{2/(1-|
1-2/p|),2} (\ell_1,\ell_p), 
\end{equation} 
and as a sort of
non-commutative analogue 
\begin{equation} 
\label{spkwapien}
\LL(\ell_1,\ui_p)=\MM_{2/{|1-2/p|},2}(\ell_1,\ui_p)=\Pi_{2/(1-|1-
2/p|),2} (\ell_1,\ui_p). 
\end{equation} 
\end{ex} 
\begin{ex}
\label{example1b} 
More generally, Pisier \cite{pisier79} showed
that for $1<p<2$   every $p$-convex and $p'$-concave Banach
lattice $X$ is a $2/{p'}$-Hilbert space, therefore
\begin{equation} 
\label{blkwapien}
\LL(\ell_1,X)=\MM_{2p/(2-p),2}(\ell_1,X)=\Pi_{p',2}(\ell_1,X), 
\end{equation} 
and if in addition $X$ is a symmetric Banach
sequence space, then by e.\,g.  \cite{arazy} it is known that
$\ui_X$ is also a $2/{p'}$-Hilbert space,  which gives
\begin{equation} 
\label{uikwapien}
\LL(\ell_1,\ui_X)=\MM_{2p/(2-p),2}(\ell_1,\ui_X)=\Pi_{p',2}(\ell_
1,\ui_X).  
\end{equation} 
\end{ex} 
\begin{ex} 
\label{matt} 
Let
$0<\theta'<\theta<1$ and $2/(2-\theta) \le q \le 2/\theta$. Then
by  \cite[Theorem~C]{matter} every $(\theta,q)$-Hilbert space
$F$ (this means that $F=[F_0,F_1]_{\theta,q}$ for some
interpolation couple $[F_0,F_1]$ and $F_1$ a Hilbert space) is
isomorphic to a subspace of a $\theta'$-Hilbert space, hence for
such a Banach space $F$ Theorem~\ref{theta1} holds for $\theta'$
instead of $\theta$. The most prominent examples for 
$(\theta,q)$-Hilbert spaces are Lorentz spaces $\ell_{p,q}$ and
their  associated unitary ideals $\ui_{p,q}$---but note that the
results for these spaces are also included in
Example~\ref{example1b}. 
\end{ex} 
\vspace{10pt} 
Trace duality
allows interesting reformulations of Theorem~\ref{theta1}.
\begin{cor} 
\label{linfty} 
Let $F$ be a $\theta$-Hilbert space,
$0 < \theta < 1$. 
\begin{enumerate}[(a)] 
\item $\LL(\ell_\infty,
F) = \MM_{\frac{2}{1-\theta},2}(\ell_\infty,F) = 
\PP_{\frac{2}{\theta},2}(\ell_\infty,F)$. 
\item Every
$\frac{2}{1-\theta}$-summing operator on $F$ factorizes through
a  Hilbert space. 
\end{enumerate} 
\end{cor} 
\proof  (b) follows
from \eqref{thetaa} by trace duality: By local techniques (see
again \cite[23.1]{df}) statement~\eqref{thetaa} in terms of
quotient ideals (see e.\,g. \cite[25.6]{df}) reads as follows:
$$ 
\PP_{\frac{2}{1-\theta}}(F, \cdot) \subset (\PP_1 \circ
\Gamma_1^{-1}) (F, \cdot),
$$ where $\Gamma_p$ for $1 \le p \le
\infty$ stands for the Banach operator  ideal of all
$T:F\rightarrow Y$ such that $F \overset{T}{\rightarrow} Y \Id
Y''$ factorizes through some $L_p(\mu)$. Hence the abstract
quotient formula from \cite[25.7]{df}  together with the trace
formula $\Pi_1^*=\mathcal{I}_\infty=\Gamma_\infty$  (see e.\,g.
\cite[6.16]{djt}) and the fact that the adjoint $\Gamma_2^*$ of 
$\Gamma_2$ is contained in $\Gamma_1 \circ \Gamma_\infty$ (a
result of Kwapie\'n, see e.\,g. \cite[7.12]{djt}) imply the
conclusion: 
$$ 
\PP_{\frac{2}{1-\theta}}(F,\cdot) \subset
(\Gamma_1 \circ \Gamma_\infty)^*  \subset \Gamma_2.
$$ 
Finally,
(a) is an immediate consequence of (b): Take $T \in
\LL(\ell_\infty, F)$ and some $S \in
\PP_{\frac{2}{1-\theta}}(F,Y) \subset \Gamma_2(F, Y)$ 
(by (b)), $Y$ some Banach space. Then $ST$ by the little Grothendieck theorem is
$2$-summing. \qed 
\par For $\theta=1$ and $F=\ell_2$ the
statements \eqref{lpkwapien} and  Corollary~\ref{linfty}~(a) are
 the ``Little Grothendieck Theorems''. The special cases 
$\LL(\ell_1,\ell_p) = \PP_{r,2}(\ell_1, \ell_p)$, $1/r = 1/2-
|1/2 - 1/p|$, and  $\LL(\ell_\infty,\ell_p)=
\PP_{p,2}(\ell_\infty,\ell_p)$, $2 \le p \le  \infty$,  are due
to Kwapie\'n \cite{kwapien} and  Lindenstrauss--Pe\l czy\'nski
\cite{lp}, respectively, whereas  Example~\ref{example1b} for
summing operators is contained in \cite{pisier79}. 
Example~\ref{matt} is due to  Lermer \cite{lermer}, and
 Corollary~\ref{linfty} in the
present form seems to be new (in the case $F=\ell_p$ see
\cite[Ex.~34.12]{df} and \cite[p.~168]{djt}). 
\begin{cor} 
Let $2
\le r,s \le \infty$ and $1/r =1/2 -1/s $. 
\vspace{-5pt}
\begin{enumerate}[(a)] 
\item $\ui_r = \PP_{r,2}(\ell_2) =
\MM_{s,2}(\ell_2)$ (isometrically). 
\item
$\PP_{r,2}(E,\ell_2)=\MM_{s,2}(E,\ell_2)$ for every Banach space
$E$. 
\end{enumerate} 
\end{cor} 
The first equality in (a) is due
to Mitiagin and was first published in \cite{kwapien} (see
e.\,g. \cite[1.d.12]{koenig} or \cite[10.3]{djt} for an
elementary proof). The second equality in (a) was proved in
\cite{CD92}---here is an alternative  proof by interpolation and
the first equality:  By Proposition~\ref{mixing2} for $\theta$ defined 
by $1/r =(1-\theta)/2$ and the first
equality the
embeddings in 
$$
\ui_r^n = [\ui_2^n,\ui_\infty^n]_\theta
=[\MM_{\infty,2}(\ell_2^n),  \MM_{2,2}(\ell_2^n)]_\theta \Id
\MM_{s,2}(\ell_2^n) \Id  \PP_{r,2}(\ell_2^n)= \ui_r^n
$$ 
all have
norm $\le 1$, and by localization this gives the claim. Now it
is easy to prove (b), which is a kind of extension of (a): Since
$\MM_{s,2} = \I_{s'} \circ \PP_2^{-1}$ (see e.\,g.
\cite[32.1]{df};  $\I_{s'}$ denotes the ideal of $s'$-integral
operators), it suffices  to show that $TS \in \I_{s'}$ whenever
$T \in \PP_{r,2}(E,\ell_2)$ and $S \in \PP_2(X,E)$. But by the 
Grothendieck--Pietsch factorization theorem (see e.\,g.
\cite[11.3]{df}) we know that $S=UV$ where $V \in \PP_2(X,H)$,
$U \in \LL(H,E)$ and  $H$ a Hilbert space. Then by (a)  and
local techniques $TU \in \PP_{r,2}(H,\ell_2) =
\MM_{s,2}(H,\ell_2)$, which gives $TS \in \I_{s'}$. 
\section[Rademacher and Gaussian $p$-summing operators on
$\ell_1$] {Rademacher and Gaussian $\boldsymbol{p}$-summing
operators on $\boldsymbol{\ell_1}$} 
For $1\le p \le 2$
Kwapie\'{n}'s result reads as follows: 
\begin{equation}
\label{asfollows}
\LL(\ell_1,\ell_p)=\Pi_{r,1}(\ell_1,\ell_p)=\Pi_{2,p}(\ell_1,\ell
_p),  \quad 1/r=3/2-1/p, 
\end{equation} 
where the last equality
follows by the usual inclusion formulas for  summing operators
(see e.\,g. \cite[10.4]{djt}). 
\par We now extend
\eqref{asfollows} within the framework of  Rademacher and
Gaussian $p$-summing operators---we will ``replace the  $2$ in
the last equality by Rademacher or Gauss functions''.  For $1\le
p \le 2$ an operator $T:X \rightarrow Y$ between  Banach spaces
$X$ and $Y$ is said to belong to the class of  {\sl Rademacher
$p$-summing operators} ($\Pi_{\mathcal{R},p}$ for short) if 
there exists  a constant $c>0$ such that for all finite
sequences $x_1, \ldots, x_n$  in $X$ 
\begin{equation}
\label{gammap} 
\left ( \int_0^1 \norm{\sum_{i=1}^n r_i(t) \cdot
Tx_i}^2 dt  \right)^{1/2} \le c \cdot \sup_{x' \in B_{X'}} \left
( \sum_{i=1}^n |\langle x',x_i  \rangle |^p \right )^{1/p},
\end{equation}  
where $(r_i)_i$ denotes as usual the sequence of
Rademacher functions;  these are all operators which transform
weakly $p$-summable sequences into almost unconditionally
summable sequences (to see this use e.\,g. \cite[12.3]{djt}).  
We write  $\pi_{\mathcal{R},p}(T)$ for the smallest constant $c$
satisfying  \eqref{gammap} and  obtain in this way  the
injective and maximal Banach operator ideal 
$(\Pi_{\mathcal{R},p},\pi_{\mathcal{R},p})$. If in
\eqref{gammap} the Rademacher functions are substituted  by a
sequence of independent Gaussian variables $(g_i)_i$, then  we
denote the resulting Banach  operator ideal of all {\em Gaussian
$p$-summing operators}  by $(\Pi_{\gamma,p},\pi_{\gamma,p})$;
for $p=2$ it was originally  introduced by \cite{linde}. Note
that $\Pi_{\mathcal{R},2}$ and $\Pi_{\gamma,2}$ coincide (see
e.\,g.  \cite[12.12]{djt}), and  if $Y$ has non-trivial cotype,
then  $\Pi_{\gamma,p}(X,Y)=\Pi_{\mathcal{R},p}(X,Y)$ for any
Banach space $X$ (use \cite[12.11 and 12.27]{djt}).  Finally,
$\LL = \Pi_{\mathcal{R},1}$ (apply the Kahane inequality
\cite[11.1]{djt}). 
\par Let $0<\theta<1$. A Banach space $F$ is
called a {\em $\theta$-type $2$  space} if there exists an
interpolation couple $[F_0,F_1]$ such that  $F_1$ has type $2$
and $F=[F_0,F_1]_\theta$. As in the case of  $\theta$-Hilbert
spaces one can always assume that $F_0 \cap F_1$ is dense  in
$F_0$ and $F_1$. 
\par Our extension of Kwapie\'{n}'s result
\eqref{asfollows} within the framework of Rademacher $p$-summing
operators now reads as follows: 
\begin{theo} 
\label{rademacher}
For $1<p<2$ and $\theta:=2/{p'}$ let $F$ be a $\theta$-type $2$
space. Then 
$$ \LL(\ell_1,F)=\Pi_{\mathcal{R},p}(\ell_1,F). 
$$
\em Note that by interpolation a $\theta$-type~$2$~space $F$,
and therefore  by $K$-convexity and duality also its dual space
$F'$, have non-trivial type,  which implies that $F$ has
non-trivial cotype. Hence, in the above theorem one may
substitute $\Pi_{\mathcal{R},p}(\ell_1,F)$ by 
$\Pi_{\gamma,p}(\ell_1,F)$.    
\end{theo} 
Before we give the
proof let us again illustrate our result by a first  example.
\begin{ex} 
Since every $\theta$-Hilbert space is a
$\theta$-type~$2$~space, $\ell_p$  and $\ui_p$ for $1<p<2$ are
$2/{p'}$-type~$2$ spaces  (see Example~\ref{example1a}). Hence,
by Theorem~\ref{rademacher} 
$$
\LL(\ell_1,\ell_p)=\Pi_{\mathcal{R},p}(\ell_1,\ell_p) \qquad
\text{and}  \qquad
\LL(\ell_1,\ui_p)=\Pi_{\mathcal{R},p}(\ell_1,\ui_p). 
$$ 
This
gives an alternative approach to Kwapie\'{n}'s result  in the
case $1<p<2$: $\ell_p$ and  $\ui_p$ both have cotype~2 (the
latter result is due to \cite{tom}),  hence, by the definition
of the Rademacher $p$-summing operators, for $V_p=\ell_p$ or
$\ui_p$  $$ \LL(\ell_1,V_p)=\Pi_{\mathcal{R},p}(\ell_1,V_p)
\subset  \Pi_{2,p}(\ell_1,V_p)=\Pi_{r,1}(\ell_1,V_p), $$  where
$1/r=3/2-1/p$ (see the preliminaries in Section~$2$ for the last
 equality). 
\end{ex} 
In our proof of Theorem~\ref{rademacher} we
interpolate between the   well-known fact  that a Banach space
$X$  has type~$2$ if and only if
$\LL(\ell_1,X)=\Pi_{\mathcal{R},2}(\ell_1,X)$  (see e.\,g.
\cite[12.10]{djt})  and the simple observation
$\LL=\Pi_{\mathcal{R},1}$ from above.  For this  we need an
analogue to  Proposition~\ref{mixing2}. We define (as an
extension of the definition of $d_\theta[\cdot,\cdot]$)  for a
finite-dimensional interpolation couple $[E_0,E_1]$,
$0<\theta<1$, and  $1 \le p_0,p_1,p_\theta \le 2$ such that
$1/{p_\theta}=(1-\theta)/{p_0} +  \theta/{p_1}$ 
$$
d_\theta[\ell_{p_0},\ell_{p_1};E_0,E_1]:= \sup_m 
\norm{\LL(\ell_{p_\theta'}^m ,[E_0,E_1]_\theta) \Id
[\LL(\ell_{p_0'}^m,E_0),\LL(\ell_{p_1'}^m,E_1)]_\theta}. 
$$ 
By
\cite{kouba} (see also \cite[Proposition~8]{dm99}) this quantity
is  always bounded, and the estimate $d_\theta[\ell_1,\ell_2;
\ell_u^n,\ell_u^n] \le 16$ for all $1 \le u \le 2$ will be
crucial for our applications. 
\begin{samepage} 
\begin{prop}
\label{binterpol} 
Let $[E_0,E_1]$ and $[F_0,F_1]$ be 
finite-dimensional interpolation couples. Then for $1 \le
p_0,p_1 \le 2$ and  $0<\theta<1$  
$$
\norm{[\Pi_{\mathcal{R},p_0}(E_0,F_0),
\Pi_{\mathcal{R},p_1}(E_1,F_1)]_\theta \Id 
\Pi_{\mathcal{R},p_\theta}([E_0,E_1]_\theta,[F_0,F_1]_\theta)}
\le  d_\theta[\ell_{p_0},\ell_{p_1};E_0,E_1], 
$$ 
where
$1/{p_\theta} = (1-\theta)/{p_0} + \theta/{p_1}$. \end{prop}
\end{samepage} 
\proof Throughout the proof we denote
$L_2:=L_2[0,1]$.  Consider for $\eta =0,1$ the bilinear mapping
$$  
\begin{array}{lccccc} 
\Phi_\eta^{m}: &
\Pi_{\mathcal{R},p_\eta}(E_\eta,F_\eta) & \times  &
\LL(\ell_{p_\eta'}^m,E_\eta) & \rightarrow & L_2(F_\eta) 
\\ & T
& \times & S & \mapsto & \sum_{i=1}^m r_i \cdot TSe_i.
\end{array} 
$$ 
By definition $\norm{\Phi_\eta^{m}} \le 1$, hence
for the  interpolated mapping  
\begin{multline*}
[\Phi_0^{m},\Phi_1^{m}]_\theta : 
[\Pi_{\mathcal{R},p_0}(E_0,F_0),
\Pi_{\mathcal{R},p_1}(E_1,F_1)]_\theta \times 
[\LL(\ell_{p_0'}^m,E_0), \LL(\ell_{p_1'}^m,E_1)]_\theta \\
\rightarrow  [L_2(F_0),L_2(F_1)]_\theta =L_2([F_0,F_1]_\theta)
\end{multline*} 
we obtain
$\norm{[\Phi_0^{m},\Phi_1^{m}]_\theta}\le 1.$ It follows that
for each $T:[E_0,E_1]_\theta \rightarrow [F_0,F_1]_\theta$  and
each $S \in \LL(\ell_{p_\theta'}^m,[E_0,E_1]_\theta)$ 
\begin{multline*} 
\bignorm{\sum_{i=1}^m r_i \cdot  TSe_i }_{
L_2([F_0,F_1]_\theta)}  
\\   \begin{split} &\le 
\norm{T}_{[\Pi_{\mathcal{R},p_0}(E_0,F_0),
\Pi_{\mathcal{R},p_1}(E_1,F_1)]_\theta}  \cdot 
\norm{S}_{[\LL(\ell_{p_0'}^m,E_0),
\LL(\ell_{p_1'}^m,E_1)]_\theta} 
\\ & \le 
d_\theta[\ell_{p_0},\ell_{p_1};E_0,E_1] \cdot 
\norm{T}_{[\Pi_{\mathcal{R},p_0}(E_0,F_0),
\Pi_{\mathcal{R},p_1}(E_1,F_1)]_\theta}  \cdot 
\norm{S}_{\LL(\ell_{p_\theta'}^m,[E_0,E_1]_\theta)}, 
\end{split}
\end{multline*} 
hence 
\begin{samepage} 
$$
\norm{T}_{\Pi_{\mathcal{R},p_\theta}([E_0,E_1]_\theta,[F_0,F_1]_\theta)} 
\le  d_\theta[\ell_{p_0},\ell_{p_1};E_0,E_1] \cdot
\norm{T}_{[\Pi_{\mathcal{R},p_0}(E_0,F_0), 
\Pi_{\mathcal{R},p_1}(E_1,F_1)]_\theta}. 
$$ 
\qed 
\end{samepage}
\par {\em Proof} of Theorem \ref{rademacher}: Let
$F=[F_0,F_1]_\theta$  where $F_1$ has type $2$ and $F_0 \cap
F_1$ is dense in $F_0$ and $F_1$  (this again implies that
$([F_0,F_1],F_0 \cap F_1,FIN(F_0 \cap F_1))$  is a cofinal
interpolation triple), and let $M \in FIN(F_0 \cap F_1)$.  An
easy application of the Kahane inequality yields $$
\norm{\LL(\ell_1^n, M_0) \Id \Pi_{\mathcal{R},1}(\ell_1^n,M_0)}
\le  \sqrt{2}, $$ and by \cite[p.~245]{djt} $$
\norm{\LL(\ell_1^n,M_1) \Id \Pi_{\mathcal{R},2}(\ell_1^n,M_1)}
\le K_G  \cdot \type{2}(F_1), $$ where $K_G$ denotes the
Grothendieck constant. Hence, by the usual  interpolation
theorem together with Kwapie\'{n}'s interpolation trick and
Proposition~\ref{binterpol}, we obtain (recall that 
$d_\theta[\ell_1,\ell_2;\ell_1^n,\ell_1^n] \le 16$) 
$$
\norm{\LL(\ell_1^n, [M_0,M_1]_\theta) \Id
\Pi_{\mathcal{R},p}(\ell_1^n, [M_0,M_1]_\theta)} \le 16 \cdot
2^{(1-\theta)/2} \cdot (K_G \cdot  \type{2}(F_1))^\theta. 
$$ 
The
claim now follows by Lemma~\ref{l1lift}. \qed 
\par For our
second example the following extrapolation theorem is crucial:
\begin{theo} 
\label{extra} 
For $1<p<2$  let $X$ be a $p$-convex
Banach function space with finite concavity. Then $X(\C)$ is a
$\theta$-type~$2$ space, $\theta:=2/{p'}$. 
\end{theo} 
For the
proof of this theorem we need some further tools. Let
$X_0(\mu),X_1(\mu)$ be Banach function spaces and $0<\theta<1$. 
 The Banach function space $X_0^{1-\theta}  X_1^\theta$ is
defined to be the  set of  functions $f \in L_0(\mu)$ for which
there exist $g \in X_0$ and $h \in X_1$ such that
$|f|=|g|^{1-\theta} \cdot |h|^\theta$, together with the norm 
$$ 
\norm{f}_{X_0^{1-\theta}  X_1^\theta}:= \inf \{ 
\norm{g}_{X_0}^{1-\theta} \cdot \norm{h}_{X_1}^\theta \, | \, 
|f|=|g|^{1-\theta} \cdot |h|^\theta, g \in X_0, h \in X_1 \}. 
$$
For $0<r<\infty$ and a Banach function space $X(\mu)$ with
$\convex{\max(1,r)}(X)=1$ we  define the  Banach function space
$(X^r,\norm{\cdot}_{X^r})$ by  
$$ X^r := \{f \in L_0(\mu) \, |
\, |f|^{1/r} \in X \} \qquad \text{and} \qquad \norm{f}_{X^r} :=
\norm{|f|^{1/r}}_X^r, \quad f \in X^r. 
$$ 
An easy calculation
(see e.\,g. \cite[Lemma~2]{defant}) shows that $X^r$ is 
$(s/r)$-convex and $(t/r)$-concave whenever  $\max(1,r) \le t,s
< \infty$ and $X$ is $s$-convex and $t$-concave. 
\par The
following crucial interpolation formula for Banach function
spaces  is due to Calder\'{o}n \cite[13.6]{calderon}: Let
$X_0(\mu),X_1(\mu)$ be  Banach function spaces such that at
least  one  is $\sigma$-order continuous. Then for each
$0<\theta<1$ isometrically 
\begin{equation} 
\label{Calderon}
[X_0(\C),X_1(\C)]_\theta = (X_0^{1-\theta}  X_1^\theta)(\C).
\end{equation} 
\par Now we are able to state the following
interpolation/extrapolation lemma.  
\begin{samepage}
\begin{lemma} 
\label{power2} 
Let $X_0(\mu)$ and $X_1(\mu)$ be
Banach function spaces, and for $0<r<1<p<\infty$  assume that
$X_0$ is $p$-convex with  $\convex{p}(X_0)=1$. Then for
$0<\theta<1$ such that $p(1-\theta)+r\theta=1$ 
$$
(X^p_0)^{1-\theta} (X_1^r)^\theta = X_0^{p(1-\theta)}
X_1^{r\theta}. 
$$ 
Moreover, if  $X$ is a $p$-convex Banach
function space with  $\convex{p}(X)=1$, then for
$\theta:=2/{p'}$ 
$$ X=(X^p)^{1-\theta}  (X^{p/2})^\theta.  
$$
\end{lemma} 
\end{samepage} 
\proof Let $V:=(X^p_0)^{1-\theta} 
(X_1^r)^\theta$ and  $W:=X_0^{p(1-\theta)}  X_1^{r\theta}$.
Then, if $f \in V$ and  $|f|=|g|^{1-\theta} \cdot |h|^\theta $
with $g \in X_0^p$ and $h \in X_1^r$, then
$|f|=(|g|^{1/p})^{p(1-\theta)} \cdot (|h|^{1/r})^{r\theta} \in
W$, and 
$$ \norm{f}_W = \norm{(|g|^{1/p})^{p(1-\theta)} \cdot
(|h|^{1/r})^{r\theta}}_W \le
\norm{|g|^{1/p}}_{X_0}^{p(1-\theta)} \cdot 
\norm{|h|^{1/r}}_{X_1}^{r\theta}= \norm{g}_{X_0^p}^{1-\theta}
\cdot \norm{h}_{X_1^r}^\theta, 
$$ 
hence $\norm{f}_W \le
\norm{f}_V$. Conversely, let $f \in W$ and 
$|f|=|g|^{p(1-\theta)} \cdot |h|^{r\theta}$ with $g \in X_0$ and
 $h \in X_1$. Then $|f|=(|g|^p)^{1-\theta} \cdot (|h|^r)^\theta
\in V$, and 
$$ \norm{f}_V =\norm{(|g|^p)^{1-\theta} \cdot
(|h|^r)^\theta}_V \le \norm{|g|^p}_{X_0^p}^{1-\theta} \cdot
\norm{|h|^r}_{X_1^r}^\theta = \norm{g}_{X_0}^{p(1-\theta)} \cdot
\norm{h}_{X_1}^{r\theta}, 
$$ 
hence $\norm{f}_V \le \norm{f}_W$.
For the rest observe that  $X=X^{1-\eta}  X^\eta$ holds
isometrically with equal norms for  each $0<\eta<1$; this
follows from the  abstract H\"older inequality (see e.\,g.
\cite[1.d.2]{lt}): Let $f \in  X^{1-\eta} X^\eta$ and
$|f|=|g|^{1-\eta} \cdot |h|^\eta$ with  $g,h \in X$. Then $|f|
=|g|^{1-\eta} \cdot |h|^\eta \in X$, and 
$$ \norm{f}_X =
\norm{|g|^{1-\eta} \cdot |h|^\eta}_X \le  \norm{g}_X^{1-\eta}
\cdot \norm{h}_X^\eta, 
$$ 
hence $\norm{f}_X \le
\norm{f}_{X^{1-\eta} X^\eta}$. Conversely, we have 
$|f|=|f|^{1-\eta} \cdot |f|^\eta \in X^{1-\eta} X^\eta$, and 
$$
\norm{f}_{X^{1-\eta}  X^\eta} \le \norm{f}_X^{1-\eta} \cdot
\norm{f}_X^\eta = \norm{f}_X. 
$$ 
Clearly $\theta:= 2/{p'}$
satisfies $p(1-\theta)+{p/2} \cdot \theta=1$.  Altogether we
obtain that 
$$ (X^p)^{1-\theta} 
(X^{p/2})^\theta=X^{p(1-\theta)}  X^{p\theta/2}=X. 
$$ 
\qed 
\par
{\em Proof} of Theorem~\ref{extra}: Without loss of generality 
we may assume that $\convex{p}(X)=1$ (see e.\,g.
\cite[1.d.8]{lt}). By assumption there exists  $p<q<\infty$ such
that $X$ is $q$-concave. Then $X^{p/2}$ is $2$-convex and 
$(2q/p)$-concave which implies that $X^{p/2}$ is $\sigma$-order
continuous (see \cite[1.a.5 and 1.a.7]{lt}) and that
$X^{p/2}(\C)$ has type~2 (see \cite[1.f.17]{lt}). Hence,
Lemma~\ref{power2} and \eqref{Calderon} imply that
$X(\C)=[X^p(\C),X^{p/2}(\C)]_\theta$ is a
$\theta$-type~$2$~space, with  $\theta:=2/{p'}$. \qed 
\begin{ex}
\label{example2} 
For $1<p<2$ let $X$ be a $p$-convex Banach
function  space with finite concavity. Then by
Theorem~\ref{extra} we know that  $X(\C)$ is a
$\theta$-type~$2$~space, $\theta:=2/{p'}$,  and  consequently
(by complexification)  
$$
\LL(\ell_1,X)=\Pi_{\mathcal{R},p}(\ell_1,X). 
$$ 
If in addition
$X$ is a symmetric Banach sequence space,  then   $\ui_X
=[\ui_{X^p}, \ui_{X^{p/2}}]_\theta$  (by \cite{arazy}) is also a
 $\theta$-type~$2$~space ($\ui_{X^{p/2}}$ has type~$2$ by 
\cite[p.~190]{gtj}), hence 
$$
\LL(\ell_1,\ui_X)=\Pi_{\mathcal{R},p}(\ell_1,\ui_X). 
$$ 
\end{ex}
\par  Recall that a Banach space $X$ has type~$2$ if and
only if $\LL(\ell_1,X)=\Pi_{\mathcal{R},2}(\ell_1,X)$; 
Example~\ref{example2} now reveals that the ideal
$\Pi_{\mathcal{R},p}$ might play a similar role for the notion
of type~$p$ ($1 <p<2$), at least for Banach function spaces or 
unitary ideals. 
\begin{samepage}  
If we  define as usual $$ p(X)
:= \sup \{ 1 \le p \le 2 \, | \, X \text{ has type } p \}, 
$$
and in addition 
$$ 
p_\mathcal{R}(X):= \sup \{ 1 \le p \le 2 \, |
\, \LL(\ell_1,X)=  \Pi_{\mathcal{R},p}(\ell_1,X) \}, 
$$ 
we
obtain the following: 
\end{samepage} 
\vspace{5pt} 
\begin{cor}
Let $X$ be a Banach function space or a unitary ideal.  Then
$p(X)=p_\mathcal{R}(X)$. 
\end{cor}  
\proof We start with the
case where $X$ is a Banach function space. Let  $p(X)>1$. Then
by  \cite[1.f.9]{lt} and \cite[1.f.13]{lt} $X$ is $p$-convex for
all  \mbox{$1 \le p < p(X)$} and $q$-concave for some $q <
\infty$, hence by  Example~\ref{example2}  the equality 
$\LL(\ell_1,X)= \Pi_{\mathcal{R},p}(\ell_1,X)$ holds for all  $1
\le p < p(X)$, and consequently $p_\mathcal{R}(X) \ge p(X)$. If 
$X=\ui_E$ with a symmetric Banach sequence space $E$ and
$p(\ui_E)>1$, then  $E$ has type~$p$ for all $1<p<p(\ui_E)$.
Consequently, as above,  $E$ is $p$-convex for all  \mbox{$1 \le
p < p(\ui_E)$} and $q$-concave for some $q < \infty$. Hence,  by
Example~\ref{example2},  the equality  $\LL(\ell_1,\ui_E)=
\Pi_{\mathcal{R},p}(\ell_1,\ui_E)$ holds for all  $1 \le p <
p(\ui_E)$, and therefore $p_\mathcal{R}(\ui_E) \ge p(\ui_E)$.
\par   Conversely, if  $p_\mathcal{R}(X)>1$, then by similar
arguments as in \cite[p.~237]{djt}  the Banach space $X$ is of
type~$p$ for all $1 \le p < p_\mathcal{R}(X)$,  hence $p(X) \ge
p_\mathcal{R}(X)$.  \qed  
\\[10pt]  The preceding result leads
to the following natural questions: Does 
$p(X)=p_\mathcal{R}(X)$ hold for every Banach space $X$? 
Furthermore: Is for $1 < p < 2$  a Banach space $X$ of type~$p$
if and only if  $\LL(\ell_1,X)=\Pi_{\mathcal{R},p}(\ell_1,X)$? 
\section[Eigenvalue distributions of nuclear operators on
$\theta$-type $2$  spaces]{Eigenvalue distributions of nuclear
operators on  $\boldsymbol{\theta}$-type $\boldsymbol{2}$
spaces} 
K\"{o}nig in \cite[p.~110]{koenig} shows that if $F$
is a Banach space  for which there exists some $2<s<\infty$
satisfying  
\begin{equation} 
\label{summing} 
\LL(\ell_1,F) =
\Pi_{s,2}(\ell_1,F), 
\end{equation} then each nuclear operator
$T$ on $F$ has $r$-th power summable  eigenvalues for some
$1<r<2$ (recall that such $T$ has always eigenvalues which are at least in
$\ell_2$, \cite[2.b.13]{koenig}).
\\[10pt] Moreover, he
conjectures that for each Banach space with non-trivial type
there exists an $s$ as in \eqref{summing}. The following theorem
gives  an affirmative answer for the subclass of all
$\theta$-type~$2$ spaces,  $0<\theta<1$.  
\begin{theo}
\label{nuclear} 
For $0<\theta<1$ let $F$ be a $\theta$-type~$2$
space. Then there exist $2<s,t<\infty$ (with $1/s+1/t=1/2$) such
that  
$$ 
\LL(\ell_1,F)=\MM_{s,2}(\ell_1,F)=\Pi_{t,2}(\ell_1,F).
$$ 
In particular, there exists $1<r<2$ such that every nuclear
operator on $F$ has $r$-th power summable eigenvalues.
\end{theo} 
\proof Let $F=[F_0,F_1]_\theta$ where $F_1$ has
type~$2$ and $F_0 \cap F_1$ is dense in $F_0$ and $F_1$,  and
let $M \in FIN(F_0 \cap F_1)$.  As in the proof of
Theorem~\ref{rademacher} we know that 
$$ \norm{\LL(\ell_1^n,M_1)
\Id \Pi_{\mathcal{R},2}(\ell_1,M_1)} \le K_G \cdot
\type{2}(F_1), 
$$ 
and since $F_1$ has type~$2$ it has cotype~$q$
for some $2 \le q<\infty$,  hence  
$$
\norm{\Pi_{\mathcal{R},2}(\ell_1^n,M_1) \Id
\Pi_{q,2}(\ell_1^n,M_1)} \le  \cotype{q}(F_1). 
$$ 
For each
$2<p<\infty$ such that $1/p +1/q>1/2$ there exists   (by the
``converse'' inclusion formula for summing/mixing operators 
\eqref{reverse}) a universal constant $C_{p,q}$ such that 
$$
\norm{\Pi_{q,2}(\ell_1^n,M_1) \Id \MM_{p,2}(\ell_1^n,M_1)} \le
C_{p,q}, 
$$ 
and altogether we obtain 
$$ \norm{\LL(\ell_1^n,M_1)
\Id \MM_{p,2}(\ell_1^n,M_1)} \le C_{p,q} \cdot K_G \cdot 
\type{2}(F_1) \cdot \cotype{q}(F_1). 
$$ 
Together with  
$$
\norm{\LL(\ell_1^n,M_0) \Id \MM_{2,2}(\ell_1^n,M_0)}=1, 
$$   
the
usual interpolation theorem, Kwapie\'{n}'s interpolation trick
and  Proposition~\ref{mixing2} this gives for $2 <s<\infty$ such
that $1/s=(1-\theta)/2 + \theta/p$  
$$
\norm{\LL(\ell_1^n,[M_0,M_1]_\theta) \Id \MM_{s,2}(\ell_1^n,
[M_0,M_1]_\theta)} \le (C_{p,q} \cdot K_G \cdot  \type{2}(F_1)
\cdot \cotype{q}(F_1))^\theta. 
$$ 
The claim now follows by
Lemma~\ref{l1lift} and the inclusion formula for summing/mixing
operators \eqref{inclusion}. \qed 
\\[15pt] Note that Banach
function spaces with non-trivial type by  Theorem~\ref{extra}
are $\theta$-type~$2$ spaces for some $\theta$. The  fact that
for a Banach function space $F$ with non-trivial type each
nuclear  operator $T$ on $F$ for some $r<2$ has $r$-th power
summable eigenvalues  is due to \cite[3.6]{pisier79}. We do not
know whether the class of all  $\theta$-type~$2$ spaces is a
proper subclass of the class of all  Banach spaces with
non-trivial type. Note that Kalton in \cite{kalton}  showed that
if $F$ has type~$p$ with $1<p<2$, then one cannot expect in
general that $F$ is a $\theta$-type~$2$ space for
$\theta:=2/{p'}$  (the ``expected'' $\theta$), but maybe it is
true for some smaller $\theta$.   
\providecommand{\bysame}{\leavevmode\hbox to3em{\hrulefill}\thinspace}

\end{document}